\documentclass[12pt]{article}
\usepackage{amssymb}
\usepackage{amsfonts}
\usepackage{amsmath}

\input{tcilatex}
\begin{document}

\begin{center}
\textbf{Caristi-Kirk type and Boyd\&Wong--Browder-Matkowski-Rus type fixed
point results in }$\mathbf{b}$\textbf{-metric spaces}

\bigskip

Radu MICULESCU and Alexandru MIHAIL

\bigskip
\end{center}

\textbf{Abstract}. {\small In this paper, based on a lemma giving a
sufficient condition for a sequence with elements from a }$b${\small -metric
space to be Cauchy, we obtain Caristi-Kirk type and
Boyd\&Wong--Browder-Matkowski-Rus type fixed point results in the framework
of }$b${\small -metric spaces. In addition, we extend Theorems 1,2 and 3
from [M. Bota,V. Ilea, E. Karapinar, O. Mle\c{s}ni\c{t}e, On }$\alpha _{\ast
}-\varphi ${\small -contractive multi-valued operators in }$b${\small %
-metric spaces and applications, Applied Mathematics \& Information
Sciences, 9 (2015), 2611-2620].}

\bigskip



\bigskip

\textbf{1.} \textbf{Introduction}

\bigskip

The notion of $b$-metric space was introduced by I. A. Bakhtin [3] and S.
Czerwik [14], [15] in connection with some problems concerning the
convergence of measurable functions with respect to measure.

In the last decades a considerable amount of fixed point results in the
framework of $b$-metric spaces were obtained (see, for example, [1], [2],
[6], [7], [8], [13], [17], [20], [21], [22], [23], [24], [25], [27], [28]
and the references therein).

In this paper we present a sufficient condition for a sequence with elements
from a $b$-metric space to be Cauchy (see Lemma 2.2). Then, using this
result, in Section 3, we obtain Caristi-Kirk type and
Boyd\&Wong--Browder-Matkowski-Rus type fixed point results in the framework
of $b$-metric spaces (see Theorems 3.1 and 3.2). Section 4 is devoted to the
presentation of a class of comparison functions, denoted by $\Gamma ^{\gamma
}$, satisfying the hypotheses of the above mentioned
Boyd\&Wong--Browder-Matkowski-Rus type fixed point result. We also point out
that $\Gamma ^{\gamma }$ is larger that Berinde's class $\Psi _{b}$. In
Section 5, we show how to extend, using Lemma 2.2, the results from [9]
concerning $\alpha _{\ast }$-$\varphi $-contractive multivalued operators by
considering comparison functions $\varphi $ from $\Gamma ^{\gamma }$ rather
than from $\Psi _{b}$.

\newpage 

\textbf{2.} \textbf{Preliminaries}

\bigskip

In this section we recall some basic facts that will be used in the sequel.

\bigskip

\textbf{Definition 2.1.} \textit{Given a nonempty set }$X$\textit{\ and a
real number }$s\in \lbrack 1,\infty )$\textit{, a function }$d:X\times
X\rightarrow \lbrack 0,\infty )$\textit{\ is called }$b$\textit{-metric if
it satisfies the following properties:}

\textit{i) }$d(x,y)=0$\textit{\ if and only if }$x=y$\textit{;}

\textit{ii) }$d(x,y)=d(y,x)$\textit{\ for all }$x,y\in X$\textit{;}

\textit{iii) }$d(x,y)\leq s(d(x,z)+d((z,y))$\textit{\ for all }$x,y,z\in X$%
\textit{.}

\textit{The pair }$(X,d)$\textit{\ is called }$b$\textit{-metric space (with
constant }$s$\textit{).}

\bigskip

Besides the classical spaces $l^{p}(\mathbb{R})$ and $L^{p}[0,1]$, where $%
p\in (0,1)$, more examples of $b$-metric spaces could be found in [2], [4],
[8], [14] and [15].

\bigskip

\textbf{Remark 2.1}. \textit{Every metric space is a }$b$\textit{-metric
space (with constant }$1$\textit{), There exist }$b$\textit{-metric spaces
which are not metric spaces (see, for example, }[2], [13] \textit{or} [21]%
\textit{).}

\bigskip

\textbf{Definition 2.2.} \textit{A sequence} $(x_{n})_{n\in \mathbb{N}}$ 
\textit{of elements from a }$b$\textit{-metric space }$(X,d)$\textit{\ is
called:}

\textit{- convergent if there exists }$l\in \mathbb{R}$\textit{\ such that }$%
\underset{n\rightarrow \infty }{\lim }d(x_{n},l)=0$;

\textit{- Cauchy if }$\underset{m,n\rightarrow \infty }{\lim }%
d(x_{m},x_{n})=0$\textit{, i.e. for every }$\varepsilon >0$\textit{\ there
exists} $n_{\varepsilon }\in \mathbb{N}$ \textit{such that }$%
d(x_{m},x_{n})<\varepsilon $\textit{\ for all} $m,n\in \mathbb{N}$, $m,n\geq
n_{\varepsilon }$.

\textit{The} $b$\textit{-metric space }$(X,d)$ \textit{is called complete if
every Cauchy sequence of elements from }$(X,d)$\textit{\ is convergent.}

\bigskip

\textbf{Remark 2.2.} \textit{A }$b$\textit{-metric space can be endowed with
the topology induced by its convergence.}

\bigskip

Using the method of mathematical induction, one can easily establish the
following result:

\bigskip

\textbf{Lemma 2.1}. \textit{For every sequence} $(x_{n})_{n\in 
\mathbb{N}
}$ \textit{of elements from a }$b$\textit{-metric space }$(X,d)$, \textit{%
with constant} $s$\textit{, the inequality }%
\begin{equation*}
d(x_{0},x_{k})\leq s^{n}\underset{i=0}{\overset{k-1}{\sum }}d(x_{i},x_{i+1})
\end{equation*}%
\textit{is valid for every }$n\in \mathbb{%
\mathbb{N}
}$\textit{\ and every }$k\in \{1,2,3,...,2^{n}-1,2^{n}\}$\textit{.}

\bigskip

The following lemma is the key ingredient in the proof of Theorems 3.1, 3.2
and 5.1.

\bigskip

\textbf{Lemma 2.2}. \textit{A sequence} $(x_{n})_{n\in 
\mathbb{N}
}$ \textit{of elements from a }$b$\textit{-metric space }$(X,d)$, \textit{%
with constant} $s$\textit{, is Cauchy provided that there exists }$\gamma
>\log _{2}s$ \textit{such that the series }$\overset{\infty }{\underset{n=1}{%
\sum }}n^{\gamma }d(x_{n},x_{n+1})$ \textit{is convergent.}

\textit{Proof}. With the notation $\alpha :=\gamma \log _{s}2>1$, we have $%
\underset{n\rightarrow \infty }{\lim }s^{(n+1)(n+2)-\alpha n^{2}}=0$ since $%
\underset{n\rightarrow \infty }{\lim }(n+1)(n+2)-\alpha n^{2}=-\infty $.
Therefore the sequence $(s^{(n+1)(n+2)-\alpha n^{2}})_{n\in 
\mathbb{N}
}$, i.e. $(\frac{s^{(n+1)(n+2)}}{(2^{n^{2}})^{\gamma }})_{n\in 
\mathbb{N}
}$, is bounded and since $\underset{x\in \mathbb{R}}{\sup }(x+1)(x+2)-\alpha
x^{2}=2+\frac{9}{4(\alpha -1)}$, we have 
\begin{equation}
\underset{n\in \mathbb{N}}{\sup }\frac{s^{(n+1)(n+2)}}{(2^{n^{2}})^{\gamma }}%
\leq s^{2+\frac{9}{4(\alpha -1)}}:=M\text{.}  \tag{1}
\end{equation}

We claim that 
\begin{equation}
d(x_{n},x_{n+m})\leq M\overset{\infty }{\underset{i=n}{\sum }}i^{\gamma
}d(x_{i},x_{i+1})\text{,}  \tag{2}
\end{equation}%
for all $m,n\in \mathbb{%
\mathbb{N}
}$.

Indeed, with the notation $l=[\sqrt{\log _{2}(m+1)}]$ (having in mind that
\linebreak $2^{l^{2}}-1\leq m<2^{(l+1)^{2}}-1$), we get%
\begin{equation*}
d(x_{n},x_{n+m})\leq sd(x_{n},x_{n+1})+sd(x_{n+1},x_{n+m})\leq
\end{equation*}%
\begin{equation*}
\leq
sd(x_{n+2^{0^{2}}-1},x_{n+2^{(0+1)^{2}}-1})+s^{2}d(x_{n+2^{1^{2}}-1},x_{n+2^{2^{2}}-1})+s^{2}d(x_{n+2^{2^{2}}-1},x_{n+m})\leq
\end{equation*}%
\begin{equation*}
...
\end{equation*}%
\begin{equation*}
\leq \overset{l-1}{\underset{i=0}{\sum }}%
s^{i+1}d(x_{n+2^{i^{2}}-1},x_{n+2^{(i+1)^{2}}-1})+s^{l}d(x_{n+2^{l^{2}}-1},x_{n+m})%
\overset{\text{Lemma 1}}{\leq }
\end{equation*}%
\begin{equation*}
\leq \overset{l-1}{\underset{i=0}{\sum }}s^{i+1}s^{(i+1)^{2}}(\underset{%
j=2^{i^{2}}}{\overset{2^{(i+1)^{2}}-1}{\sum }}%
d(x_{n+j-1},x_{n+j}))+s^{l+1}s^{(l+1)^{2}}\underset{j=2^{l^{2}}}{\overset{%
2^{(l+1)^{2}}-1}{\sum }}d(x_{n+j-1},x_{n+j}))=
\end{equation*}%
\begin{equation*}
=\overset{l}{\underset{i=0}{\sum }}s^{i+1}s^{(i+1)^{2}}(\underset{j=2^{i^{2}}%
}{\overset{2^{(i+1)^{2}}-1}{\sum }}d(x_{n+j-1},x_{n+j}))\leq
\end{equation*}%
\begin{equation*}
\leq \overset{\infty }{\underset{i=0}{\sum }}\frac{s^{(i+1)(i+2)}}{%
(2^{i^{2}})^{\gamma }}(\underset{j=2^{i^{2}}-1}{\overset{2^{(i+1)^{2}}-1}{%
\sum }}(j+1)^{\gamma }d(x_{n+j-1},x_{n+j}))\overset{(1)}{\leq }
\end{equation*}%
\begin{equation*}
\leq M\overset{\infty }{\underset{i=0}{\sum }}(i+1)^{\gamma
}d(x_{n+i},x_{n+i+1})\leq M\overset{\infty }{\underset{i=0}{\sum }}%
(n+i)^{\gamma }d(x_{n+i},x_{n+i+1})=
\end{equation*}%
\begin{equation*}
=M\overset{\infty }{\underset{i=n}{\sum }}i^{\gamma }d(x_{i},x_{i+1})\text{.}
\end{equation*}

As the series $\overset{\infty }{\underset{n=1}{\sum }}n^{\gamma
}d(x_{n},x_{n+1})$ is convergent, from $(2)$, we infer that $(x_{n})_{n\in 
\mathbb{N}
}$ is Cauchy. $\square $

\bigskip

Using the comparison test, from the above lemma, we obtain the following two
results.

\bigskip

\textbf{Corollary 2.1}. \textit{A sequence} $(x_{n})_{n\in 
\mathbb{N}
}$ \textit{of elements from a }$b$\textit{-metric space }$(X,d)$\textit{,
with constant} $s$\textit{, is Cauchy provided that there exist }$\gamma
>\log _{2}s$ \textit{and a sequence }$(a_{n})_{n\in 
\mathbb{N}
}$ \textit{of positive real numbers such that:}

\textit{a)} \textit{the series }$\overset{\infty }{\underset{n=1}{\sum }}%
a_{n}d(x_{n},x_{n+1})$ \textit{is convergent;}

\textit{b) }$\underline{lim}\frac{a_{n}}{n^{\gamma }}>0$.

\bigskip

\textbf{Corollary 2.2}. \textit{A sequence} $(x_{n})_{n\in 
\mathbb{N}
}$ \textit{of elements from a }$b$\textit{-metric space }$(X,d)$\textit{\ is
Cauchy provided that there exists }$\alpha >1$\textit{\ such that} \textit{%
the series}\linebreak \textit{\ }$\overset{\infty }{\underset{n=1}{\sum }}%
\alpha ^{n}d(x_{n},x_{n+1})$ \textit{is convergent.}

\bigskip

\textbf{3. Caristi-Kirk type and Boyd\&Wong--Browder-Matkowski-Rus type
fixed point results in }$\mathbf{b}$\textbf{-metric spaces}

\bigskip

In this section, using Lemma 2.2, we obtain two fixed point theorems in the
framework of $b$-metric spaces.

\bigskip

Our first result is a Caristi-Kirk type fixed point result (see [12], [16]
and [18]).

\bigskip

\textbf{Theorem 3.1}. \textit{Let} $(X,d)$\textit{\ be a complete }$b$%
\textit{-metric space, }$\varphi :X\rightarrow \lbrack 0,\infty )$, $%
f:X\rightarrow X$ \textit{and} $\alpha >1$ \textit{such that:}

\textit{a) }$f$\textit{\ is continuous;}

\textit{b) }$d(x,f(x))\leq \varphi (x)-\alpha \varphi (f(x))$\textit{\ for
every }$x\in X$\textit{.}

\textit{Then, for every }$x_{0}\in X$\textit{, the sequence }$%
(f^{[n]}(x_{0}))_{n\in \mathbb{N}}$\textit{\ is convergent and its limit is
a fixed point of }$f$\textit{.}

\textit{Proof}. With the notation $x_{n}:=f^{[n]}(x_{0})$, according to b),
we have $d(x_{n},x_{n+1})\leq \varphi (x_{n})-\alpha \varphi (f(x_{n+1}))$,
so%
\begin{equation*}
\alpha ^{n}d(x_{n},x_{n+1})\leq \alpha ^{n}\varphi (x_{n})-\alpha
^{n+1}\varphi (f(x_{n+1}))\text{,}
\end{equation*}%
for every $n\in \mathbb{N}$ and consequently the series $\overset{\infty }{%
\underset{n=1}{\sum }}\alpha ^{n}d(x_{n},x_{n+1})$ is convergent, its
partial sums being between $0$ and $\alpha \varphi (x_{1})$. According to
Corollary 2.2, the sequence $(x_{n})_{n\in \mathbb{N}}$ is convergent and if
we denote by $u$ its limit, then passing to limit as $n\rightarrow \infty $
in the relation $x_{n+1}=f(x_{n})$, based on a), we infer that $f(u)=u$,
i.e. $u$ is a fixed point of $f$. $\square $

\bigskip

\textbf{Remark 3.1}. \textit{The above result gives us a sufficient
condition for }$f$\textit{\ to be a weakly Picard operator.}

\bigskip

Our second result is a Boyd\&Wong--Brower-Matkowski-Rus type fixed point
result (see [10], [11], [19] and [26]).

\bigskip

\textbf{Theorem 3.2}. \textit{Let} $(X,d)$\textit{\ be a complete }$b$%
\textit{-metric space, with constant }$s$\textit{, }$\gamma >\log _{2}s$,%
\textit{\ }$\varphi :[0,\infty )\rightarrow \lbrack 0,\infty )$ \textit{and} 
$f:X\rightarrow X$ \textit{such that:}

\textit{a) }$\varphi (r)<r$\textit{\ for every }$r>0$\textit{;}

\textit{b) the series} $\overset{\infty }{\underset{n=1}{\sum }}n^{\gamma
}\varphi ^{\lbrack n]}(r)$ \textit{is convergent for every }$r>0$\textit{;}

\textit{c) }$f$\textit{\ is a }$\varphi $\textit{-contraction, i.e. }$%
d(f(x),f(y))\leq \varphi (d(x,y))$\textit{\ for all }$x,y\in X$\textit{.}

\textit{Then }$f$ \textit{has a unique fixed point }$u$\textit{\ and the
sequence }$(f^{[n]}(x_{0}))_{n\in \mathbb{N}}$\textit{\ is convergent to }$u$
\textit{for every }$x_{0}\in X$\textit{.}

\textit{Proof}. With the notation $x_{n}:=f^{[n]}(x_{0})$, taking into
account c), we have%
\begin{equation}
n^{\gamma }d(x_{n},x_{n+1})\leq n^{\gamma }\varphi ^{\lbrack
n]}(d(x_{1},x_{0}))\text{,}  \tag{1}
\end{equation}%
for every $n\in \mathbb{N}$.

Based on $(1)$, b) and the comparison test, we conclude that the series $%
\overset{\infty }{\underset{n=1}{\sum }}n^{\gamma }d(x_{n},x_{n+1})$ is
convergent and Lemma 2.2 assures us that the sequence $(x_{n})_{n\in \mathbb{%
N}}$ is convergent. If we denote by $u$ its limit, then passing to limit as $%
n\rightarrow \infty $ in the relation $x_{n+1}=f(x_{n})$, since $f$ is
continuos (see c)), we infer that $f(u)=u$, i.e. $u$ is a fixed point of $f$.

Now let us prove that $u$ is unique.

Indeed, if there exists $v\in X\smallsetminus \{u\}$ having the property
that $f(v)=v$, then we get the following contradiction: $d(u,v)=d(f(u),f(v))%
\leq \varphi (d(u,v))\overset{a)}{<}d(u,v)$. $\square $

\bigskip

\textbf{Remark 3.2}. We have%
\begin{equation*}
d(x_{n},u)\leq sd(x_{n},x_{n+m})+sd(x_{n+m},u)\overset{(2)\text{ from the
proof of Lemma 2.2}}{\leq }
\end{equation*}%
\begin{equation*}
\leq sM\overset{\infty }{\underset{i=n}{\sum }}i^{\gamma
}d(x_{i},x_{i+1})+sd(x_{n+m},u)\text{,}
\end{equation*}%
for all $m,n\in \mathbb{N}$. By passing to limit as $m\rightarrow \infty $
in the above inequality we get the following \textit{estimation of the speed
of convergence for} $(f^{[n]}(x_{0}))_{n\in \mathbb{N}}$:%
\begin{equation*}
d(x_{n},u)\leq sM\overset{\infty }{\underset{i=n}{\sum }}i^{\gamma
}d(x_{i},x_{i+1})\text{,}
\end{equation*}%
for every $n\in \mathbb{N}$.

\bigskip

\textbf{Remark 3.3}. \textit{The above result gives us a sufficient
condition for }$f$\textit{\ to be a Picard operator.}

\bigskip

\textbf{4.} \textbf{Some classes of comparison functions satisfying
conditions a) and b) from the hypotheses of Theorem 3.2}

\bigskip

In this section we introduce and study the class $\Gamma ^{\gamma }$ of $%
x^{\gamma }$-summable comparison functions - which is larger that Berinde's
class $\Psi _{b}$ - and whose elements satisfy conditions a) and b) from the
hypotheses of Theorem 3.2.

\bigskip

\textbf{Definition 4.1}. \textit{A function} $\varphi :[0,\infty
)\rightarrow \lbrack 0,\infty )$ \textit{is called a comparison function if:}

\textit{i) }$\varphi $\textit{\ is increasing;}

\textit{ii) }$\underset{n\rightarrow \infty }{\lim }\varphi ^{\lbrack
n]}(r)=0$\textit{\ for every }$r\in \lbrack 0,\infty )$\textit{.}

\bigskip

\textbf{Remark 4.1}. \textit{Every comparison function }$\varphi $\textit{\
has the property that}\linebreak \textit{\ }$\varphi (r)<r$\textit{\ for
every }$r\in (0,\infty )$\textit{.}

\bigskip

\textbf{Definition 4.2}. \textit{A function} $\varphi :[0,\infty
)\rightarrow \lbrack 0,\infty )$ \textit{is called a }$x^{\gamma }$\textit{%
-summable comparison function, where }$\gamma >0$\textit{, if:}

\textit{i) }$\varphi $\textit{\ is increasing;}

\textit{ii) the series} $\overset{\infty }{\underset{n=1}{\sum }}n^{\gamma
}\varphi ^{\lbrack n]}(r)$ \textit{is convergent for every }$r\in \lbrack
0,\infty )$\textit{.}

\textit{We denote the family of} $x^{\gamma }$\textit{-summable comparison
functions by }$\Gamma ^{\gamma }$.

\bigskip

\textbf{Remark 4.2}. \textit{Every }$\varphi \in \Gamma ^{\gamma }$\textit{,
where }$\gamma >0$\textit{, is a comparison function, so it satisfies
conditions a) and b) from the hypotheses of Theorem 3.2. Hence we are
interested in finding sufficient conditions for a comparison function }$%
\varphi $\textit{\ to belong to }$\Gamma ^{\gamma }$\textit{.}

\bigskip

\textbf{Definition 4.3}. \textit{Given }$\alpha >1$\textit{, we denote by }$%
\Gamma _{\alpha }$ \textit{the family of all comparison functions} $\varphi $
\textit{for which there exist }$a>0$\textit{\ and }$\varepsilon >0$\textit{\
such that }$\varphi (x)\leq x-ax^{\alpha }$\textit{\ for every }$x\in
\lbrack 0,\varepsilon ]$\textit{.}

\bigskip

\textbf{Lemma 4.1}. \textit{Let us consider} $\alpha >1$, $\varepsilon >0$ 
\textit{and }$a>0$\textit{\ such that }$x-ax^{\alpha }\geq 0$\textit{\ for
every }$x\in \lbrack 0,\varepsilon ]$\textit{.} \textit{Then the sequence} $%
(x_{n})_{n\in \mathbb{N}}$, \textit{given by }$x_{0}\in \lbrack
0,\varepsilon ]$\textit{\ and }$x_{n+1}=x_{n}-ax_{n}^{\alpha }$\textit{\ for
every }$n\in \mathbb{N}$\textit{, has the following property:} 
\begin{equation*}
\underset{n\rightarrow \infty }{\lim }\frac{x_{n}}{(\frac{1}{n})^{\frac{1}{%
\alpha -1}}}=(\frac{1}{a(\alpha -1)})^{\frac{1}{\alpha -1}}\text{.}
\end{equation*}

\textit{Proof}. It is clear that $(x_{n})_{n\in \mathbb{N}}$ is decreasing
and $\underset{n\rightarrow \infty }{\lim }x_{n}=0$, so $(x_{n}^{1-\alpha
})_{n\in \mathbb{N}}$ is increasing and $\underset{n\rightarrow \infty }{%
\lim }x_{n}^{1-\alpha }=\infty $.

As 
\begin{equation*}
\underset{n\rightarrow \infty }{\lim }\frac{(n+1)-n}{x_{n+1}^{1-\alpha
}-x_{n}^{1-\alpha }}=\underset{n\rightarrow \infty }{\lim }\frac{1}{%
x_{n}^{1-\alpha }((\frac{x_{n+1}}{x_{n}})^{1-\alpha }-1)}=
\end{equation*}%
\begin{equation*}
=\underset{n\rightarrow \infty }{\lim }\frac{1}{x_{n}^{1-\alpha
}((1-ax_{n}^{\alpha -1})^{1-\alpha }-1)}=\underset{n\rightarrow \infty }{%
\lim }\frac{1}{-a\frac{(1-ax_{n}^{\alpha -1})^{1-\alpha }-1}{-ax_{n}^{\alpha
-1}}}=\frac{1}{a(\alpha -1)}\text{,}
\end{equation*}%
in virtue of Stolz-Cesaro lemma we obtain that $\underset{n\rightarrow
\infty }{\lim }\frac{n}{x_{n}^{1-\alpha }}=\frac{1}{a(\alpha -1)}$,
i.e.\linebreak\ $\underset{n\rightarrow \infty }{\lim }\frac{x_{n}}{(\frac{1%
}{n})^{\frac{1}{\alpha -1}}}=(\frac{1}{a(\alpha -1)})^{\frac{1}{\alpha -1}}$%
. $\square $

\bigskip

\textbf{Lemma 4.2}. \textit{Given} $\alpha >1$\textit{, for every} $\varphi
\in \Gamma _{\alpha }$\textit{\ and every }$r\geq 0$, \textit{we have}%
\begin{equation*}
\overline{\lim }\frac{\varphi ^{\lbrack n]}(r)}{(\frac{1}{n})^{\frac{1}{%
\alpha -1}}}\in \lbrack 0,\infty )\text{.}
\end{equation*}

\textit{Proof}. Taking into account ii) from the definition of a comparison
function, there exists $n_{0}\in \mathbb{N}$ such that $\varphi ^{\lbrack
n_{0}]}(r)\in \lbrack 0,\varepsilon ]$. Since $\varphi \in \Gamma _{\alpha }$%
, we infer that $\varphi ^{\lbrack n+n_{0}]}(r)\leq x_{n}$ for every $n\in 
\mathbb{N}$, where $(x_{n})_{n\in \mathbb{N}}$ is given by $x_{0}=\varphi
^{\lbrack n_{0}]}(r)\in \lbrack 0,\varepsilon ]$ and $x_{n+1}=x_{n}-ax_{n}^{%
\alpha }$ for every\textit{\ }$n\in \mathbb{N}$\textit{.} Consequently%
\textit{\ }%
\begin{equation*}
\overline{\lim }\frac{\varphi ^{\lbrack n+n_{0}]}(r)}{(\frac{1}{n+n_{0}})^{%
\frac{1}{\alpha -1}}}=\overline{\lim }\frac{\varphi ^{\lbrack n+n_{0}]}(r)}{(%
\frac{1}{n})^{\frac{1}{\alpha -1}}}\frac{1}{(\frac{n}{n+n_{0}})^{\frac{1}{%
\alpha -1}}}\leq
\end{equation*}%
\begin{equation*}
\leq \underset{n\rightarrow \infty }{\lim }\frac{x_{n}}{(\frac{1}{n})^{\frac{%
1}{\alpha -1}}}(\frac{n+n_{0}}{n})^{\frac{1}{\alpha -1}}\overset{\text{Lemma
4.1}}{=}(\frac{1}{a(\alpha -1)})^{\frac{1}{\alpha -1}}\text{,}
\end{equation*}%
so $\overline{\lim }\frac{\varphi ^{\lbrack n]}(r)}{(\frac{1}{n})^{\frac{1}{%
\alpha -1}}}\in \lbrack 0,\infty )$. $\square $

\bigskip

Our next result provides a sufficient condition for the validity of the
inclusion $\Gamma _{\alpha }\subseteq \Gamma ^{\gamma }$.

\bigskip

\textbf{Proposition 4.1}. \textit{For every }$\alpha \in (1,2)$ \textit{and} 
$\gamma \in (0,\frac{2-\alpha }{\alpha -1})$\textit{, we have}\linebreak\ $%
\Gamma _{\alpha }\subseteq \Gamma ^{\gamma }$.

\textit{Proof}. It suffices to prove that the series\textit{\ }$\overset{%
\infty }{\underset{n=1}{\sum }}n^{\gamma }\varphi ^{\lbrack n]}(r)$ is
convergent for every $r\geq 0$.

In virtue of Lemma 4.2, there exists $n_{0}\in \mathbb{N}$ such that $\frac{%
\varphi ^{\lbrack n]}(r)}{(\frac{1}{n})^{\frac{1}{\alpha -1}}}\leq C:=%
\overline{\lim }\frac{\varphi ^{\lbrack n]}(r)}{(\frac{1}{n})^{\frac{1}{%
\alpha -1}}}+1\in \mathbb{R}$, i.e. $n^{\gamma }\varphi ^{\lbrack
n]}(r)=n^{\gamma }(\frac{1}{n})^{\frac{1}{\alpha -1}}\frac{\varphi ^{\lbrack
n]}(r)}{(\frac{1}{n})^{\frac{1}{\alpha -1}}}\leq C\frac{1}{n^{\frac{1}{%
\alpha -1}-\gamma }}$ for every $n\in \mathbb{N}$, $n\geq n_{0}$. Since the
series $\overset{\infty }{\underset{n=1}{\sum }}\frac{1}{n^{\frac{1}{\alpha
-1}-\gamma }}$ is convergent (as $\frac{1}{\alpha -1}-\gamma >1$), based on
the comparison test, we conclude that the series\textit{\ }$\overset{\infty }%
{\underset{n=1}{\sum }}n^{\gamma }\varphi ^{\lbrack n]}(r)$ is convergent. $%
\square $

\bigskip

Now we provide some other sufficient conditions for a comparison function $%
\varphi $ to belong to $\Gamma ^{\gamma }$.

\bigskip

Let us suppose that for the comparison function $\varphi $ there exist the
sequences $(a_{n})_{n\in \mathbb{N}}$ and $(b_{n})_{n\in \mathbb{N}}$ such
that:

i) $a_{n}\in (0,1)$ and $b_{n}\in (0,\infty )$ for every $n\in \mathbb{N}$;

ii) $\varphi ^{\lbrack n+1]}(r)\leq a_{n}\varphi ^{\lbrack n]}(r)+b_{n}$ for
every $n\in \mathbb{N}$ and every $r\geq 0$.

Then 
\begin{equation*}
\varphi ^{\lbrack n+1]}(r)\leq
a_{n}a_{n-1}...a_{2}a_{1}r+a_{n}a_{n-1}...a_{2}b_{1}+a_{n}a_{n-1}...a_{3}b_{2}+...+a_{n}a_{n-1}b_{n-2}+a_{n}b_{n-1}+b_{n}%
\text{,}
\end{equation*}%
so, with the convention that $b_{0}=r$ and $\underset{j=n+1}{\overset{n}{%
\dprod }}a_{j}=1$, we have 
\begin{equation*}
n^{\gamma }\varphi ^{\lbrack n+1]}(r)\leq n^{\gamma }\underset{k=0}{\overset{%
n}{\sum }}b_{k}\underset{j=k+1}{\overset{n}{\dprod }}a_{j}\text{,}
\end{equation*}%
for every $n\in \mathbb{N}$, $r\geq 0$ and $\gamma >0$.

Consequently, \textit{a sufficient condition for }$\varphi $\textit{\ to
belong to }$\Gamma ^{\gamma }$\textit{\ is the convergence of the series }$%
\underset{n=0}{\overset{\infty }{\sum }}(n^{\gamma }\underset{k=0}{\overset{n%
}{\sum }}b_{k}\underset{j=k+1}{\overset{n}{\dprod }}a_{j})$\textit{, i.e. of
the series }$\underset{k=0}{\overset{\infty }{\sum }}(b_{k}\underset{n=k}{%
\overset{\infty }{\sum }}n^{\gamma }\underset{j=k+1}{\overset{n}{\dprod }}%
a_{j})$\textit{.}

\bigskip

Now we are going to take a closer look on this sufficient condition in two
particular cases.

\bigskip

The first particular case is the one for which the sequence $(a_{n})_{n\in 
\mathbb{N}}$ is constant (so there exists $a\in (0,1)$ such that $a_{n}=a$
for every $n\in \mathbb{N}$).

We claim that, in this case, the series $\underset{k=0}{\overset{\infty }{%
\sum }}(b_{k}\underset{n=k}{\overset{\infty }{\sum }}n^{\gamma }a^{n-k})$ is
convergent if and only if the series $\underset{k=0}{\overset{\infty }{\sum }%
}k^{\gamma }b_{k}$ is convergent.

The implication "$\Rightarrow $" is obvious as $k^{\gamma }\leq \underset{n=k%
}{\overset{\infty }{\sum }}n^{\gamma }a^{n-k}$ for every $k\in \mathbb{N}$.

For the implication "$\Leftarrow $" let us note that $\underset{k=0}{\overset%
{\infty }{\sum }}(b_{k}\underset{n=k}{\overset{\infty }{\sum }}n^{\gamma
}a^{n-k})$ can be rewritten as $\underset{k=0}{\overset{\infty }{\sum }}%
(b_{k}\underset{j=0}{\overset{\infty }{\sum }}(k+j)^{\gamma }a^{j})$ and
that there exists $C_{\gamma }\in \mathbb{R}$ such that $(k+j)^{\gamma }\leq
C_{\gamma }(k^{\gamma }+j^{\gamma })$ for every $j,k\in \mathbb{N}$. As the
series $\underset{k=0}{\overset{\infty }{\sum }}k^{\gamma }b_{k}$ is
convergent, in virtue of the comparison test, we infer that the series $%
\underset{k=0}{\overset{\infty }{\sum }}(b_{k}\underset{j=0}{\overset{\infty 
}{\sum }}C_{\gamma }k^{\gamma }a^{j})$ and $\underset{k=0}{\overset{\infty }{%
\sum }}(b_{k}\underset{j=0}{\overset{\infty }{\sum }}C_{\gamma }j^{\gamma
}a^{j})$ are convergent (as, with the notation $C_{a,\gamma }:=\underset{j=0}%
{\overset{\infty }{\sum }}j^{\gamma }a^{j}$, we have $b_{k}\underset{j=0}{%
\overset{\infty }{\sum }}C_{\gamma }k^{\gamma }a^{j}\leq b_{k}k^{\gamma }%
\frac{C_{\gamma }}{1-a}$ and $b_{k}\underset{j=0}{\overset{\infty }{\sum }}%
C_{\gamma }j^{\gamma }a^{j}\leq b_{k}k^{\gamma }C_{\gamma }C_{a,\gamma }$
for every $k\in \mathbb{N}$). Consequently $\underset{k=0}{\overset{\infty }{%
\sum }}(b_{k}\underset{j=0}{\overset{\infty }{\sum }}C_{\gamma }(k^{\gamma
}+j^{\gamma })a^{j})$ is convergent, and, based on the comparison test, we
conclude that $\underset{k=0}{\overset{\infty }{\sum }}(b_{k}\underset{j=0}{%
\overset{\infty }{\sum }}(k+j)^{\gamma }a^{j})$, i.e. $\underset{k=0}{%
\overset{\infty }{\sum }}(b_{k}\underset{n=k}{\overset{\infty }{\sum }}%
n^{\gamma }a^{n-k})$, is convergent. So we proved the following:

\bigskip

\textbf{Proposition 4.2}. \textit{A comparison function }$\varphi $\textit{\
for which there exist}\linebreak \textit{\ }$a\in (0,1)$\textit{\ and }$%
b_{n}\in (0,\infty )$\textit{, }$n\in \mathbb{N}$\textit{, such that }$%
\varphi ^{\lbrack n+1]}(r)\leq a\varphi ^{\lbrack n]}(r)+b_{n}$\textit{\ for
every }$n\in \mathbb{N}$\textit{\ and every }$r\geq 0$\textit{, belongs to }$%
\Gamma ^{\gamma }$\textit{, where }$\gamma >0$\textit{, provided that the
series }$\underset{n=0}{\overset{\infty }{\sum }}n^{\gamma }b_{n}$\textit{\
is convergent}.

\bigskip

The second particular case is the one for which there exists an increasing
sequence $(c_{n})_{n\in \mathbb{N}}$ converging to $\infty $ such that $%
a_{n}=\frac{c_{n-1}}{c_{n}}$ for every $n\in \mathbb{N}$.

A sufficient condition for the convergence of the series $\underset{k=0}{%
\overset{\infty }{\sum }}(b_{k}\underset{n=k}{\overset{\infty }{\sum }}%
n^{\gamma }\frac{c_{k}}{c_{n}})$, i.e. $\underset{k=0}{\overset{\infty }{%
\sum }}(b_{k}c_{k}\underset{n=k}{\overset{\infty }{\sum }}\frac{n^{\gamma }}{%
c_{n}})$, is the convergence of the series $\underset{n=0}{\overset{\infty }{%
\sum }}\frac{n^{\gamma }}{c_{n}}$ and $\underset{k=0}{\overset{\infty }{\sum 
}}b_{k}c_{k}$. This happens, for example, if $c_{n}=n^{\varepsilon +1+\gamma
}$ for every $n\in \mathbb{N}$, where $\varepsilon >0$ and the series $%
\underset{k=0}{\overset{\infty }{\sum }}b_{k}k^{\varepsilon +1+\gamma }$ is
convergent. So we proved the following:

\bigskip

\textbf{Proposition 4.3}. \textit{A comparison function }$\varphi $\textit{\
for which there exist }$\varepsilon >0$\textit{\ and }$b_{n}\in (0,\infty )$%
\textit{, }$n\in \mathbb{N}$\textit{, such that }$\varphi ^{\lbrack
n+1]}(r)\leq (\frac{n-1}{n})^{\varepsilon +1+\gamma }\varphi ^{\lbrack
n]}(r)+b_{n}$\textit{\ for every }$n\in \mathbb{N}$\textit{\ and every }$%
r\geq 0$\textit{, belongs to }$\Gamma ^{\gamma }$\textit{, where }$\gamma >0$%
, \textit{provided that the series }$\underset{n=0}{\overset{\infty }{\sum }}%
b_{n}n^{\varepsilon +1+\gamma }$\textit{\ is convergent}.

\bigskip

Let us recall (see [5]) the following:

\bigskip

\textbf{Definition 4.4}. \textit{For a given} $b>1$, \textit{by }$\Psi _{b}$%
\textit{\ we understand the class of functions }$\varphi :[0,\infty
)\rightarrow \lbrack 0,\infty )$\textit{\ such that:}

\textit{i) }$\varphi $\textit{\ is increasing;}

\textit{ii) there exist }$a\in (0,1)$\textit{\ and a convergent series }$%
\underset{n=0}{\overset{\infty }{\sum }}b_{n}$\textit{, where }$b_{n}\in
(0,\infty )$\textit{\ for every} $n\in \mathbb{N}$\textit{, such that} $%
b^{n+1}\varphi ^{\lbrack n+1]}(r)\leq ab^{n}\varphi ^{\lbrack n]}(r)+b_{n}$ 
\textit{for every }$n\in \mathbb{N}$\textit{\ and every }$r\geq 0$.

\bigskip

\textbf{Proposition 4.4.}\textit{\ }$\Psi _{b}\subseteq \Gamma ^{\gamma }$%
\textit{\ for every }$\gamma >1$\textit{.}

\textit{Proof}. If $\varphi \in \Psi _{b}$, then, with the notation $x_{n}=%
\frac{b_{n}}{b^{n+1}}$ for every $n\in \mathbb{N}$, we have $\varphi
^{\lbrack n+1]}(r)\leq \frac{a}{b}\varphi ^{\lbrack n]}(r)+x_{n}$ for every $%
n\in \mathbb{N}$ and the series $\underset{n=0}{\overset{\infty }{\sum }}%
n^{\gamma }x_{n}$, i.e. $\underset{n=0}{\overset{\infty }{\sum }}\frac{%
n^{\gamma }}{b^{n+1}}b_{n}$\textit{, }is convergent for every $\gamma >0$.
Then, according to Proposition 4.2, we conclude that $\varphi \in \Gamma
^{\gamma }$. $\square $

\bigskip

\textbf{An example}. We present an example of comparison function $\varphi $
which belongs to $\Gamma ^{\gamma }$ for every $\gamma \in (0,2)$, but not
to $\Psi _{b}$, for every $b>1$.

Let us consider the comparison function $\varphi :[0,\infty )\rightarrow
\lbrack 0,\infty )$, given by 
\begin{equation*}
\varphi (x)=\{%
\begin{array}{cc}
x-x^{\frac{4}{3}}\text{,} & x\in \lbrack 0,\frac{27}{64}] \\ 
\frac{27}{256}\text{,} & x\in (\frac{27}{64},\infty )%
\end{array}%
\text{.}
\end{equation*}

\bigskip

\textbf{Claim 1}. $\varphi \in \Gamma ^{\gamma }$ \textit{for every} $\gamma
\in (0,2)$.

\textit{Justification of claim 1}. Note that $\varphi \in \Gamma _{\frac{4}{3%
}}$ (just take $a=1$ and $\varepsilon =\frac{27}{64}$ in Definition 4.3),
so, according to Proposition 4.1, $\varphi \in \Gamma ^{\gamma }$ for every $%
\gamma \in (0,2)$.

\bigskip

\textbf{Claim 2}. \textit{There is no} $b>1$ \textit{such that} $\varphi \in
\Psi _{b}$.

\textit{Justification of claim 2}. If this is not the case, there exist $%
a\in (0,1)$\ and a convergent series $\underset{n=0}{\overset{\infty }{\sum }%
}b_{n}$, where $b_{n}\in (0,\infty )$\ for every $n\in \mathbb{N}$, such
that $b^{n+1}\varphi ^{\lbrack n+1]}(r)\leq ab^{n}\varphi ^{\lbrack
n]}(r)+b_{n}$ for every $n\in \mathbb{N}$\ and every $r\geq 0$. Hence, for a
fixed $r_{0}\in (0,\frac{27}{64})$, we have $b^{n+1}\varphi ^{\lbrack
n+1]}(r_{0})\leq ab^{n}\varphi ^{\lbrack n]}(r_{0})+b_{n}$, and therefore%
\newline
$\frac{b^{n+1}}{n^{3}}[(\frac{n}{n+1})^{3}(n+1)^{3}\varphi ^{\lbrack
n+1]}(r_{0})-\frac{a}{b}n^{3}\varphi ^{\lbrack n]}(r_{0})]\leq b_{n}$ for
every $n\in \mathbb{N}$. Since $\underset{n\rightarrow \infty }{\lim }%
n^{3}\varphi ^{\lbrack n]}(r_{0})=27$ (see Lemma 4.1), by passing to limit
as $n\rightarrow \infty $ in the previous inequality, we get that $\underset{%
n\rightarrow \infty }{\lim }b_{n}=\infty $. This is a contradiction since
the series $\underset{n=0}{\overset{\infty }{\sum }}b_{n}$ is convergent.

\bigskip

\textbf{5}. \textbf{A fixed point theorem for }$\alpha _{\ast }$-$\varphi $%
\textbf{-contractive multivalued operators in }$\mathbf{b}$\textbf{-metric
spaces}

\bigskip

In this section, inspired by the ideas from [9], we present a fixed point
theorem for $\alpha _{\ast }$-$\varphi $-contractive multivalued operator in 
$b$-metric spaces.

\bigskip

First of all let us recall some notions.

\bigskip

\textbf{Definition 5.1. }\textit{For} $T:X\rightarrow \mathcal{P}%
(X):=\{Y\mid Y\subseteq X\}$ \textit{and}\linebreak\ $\alpha :X\times
X\rightarrow \lbrack 0,\infty )$\textit{, where }$(X,d)$\textit{\ is a }$b$%
\textit{-metric space, we say that }$T$ \textit{is }$\alpha _{\ast }$\textit{%
-admissible if, for all }$x,y\in X$\textit{, the following implication is
valid: }$\alpha (x,y)\geq 1\Rightarrow \alpha _{\ast }(T(x),T(y)):=\inf
\{\alpha (u,v)\mid u\in T(x),v\in T(y)\}\geq 1$\textit{.}

\bigskip

\textbf{Definition 5.2. }\textit{For} $T:X\rightarrow \mathcal{P}%
_{cl}(X):=\{Y\in \mathcal{P}(X)\mid Y$ \textit{is closed}$\}$, $\alpha
:X\times X\rightarrow \lbrack 0,\infty )$ \textit{and }$\varphi \in \Gamma
^{\gamma }$\textit{, where }$(X,d)$\textit{\ is a }$b$\textit{-metric space
and }$\gamma >1$\textit{, we say that }$T$ \textit{is an }$\alpha _{\ast }$-$%
\varphi $\textit{-contractive multivalued operator of type }$(b)$ \textit{if 
}$\alpha _{\ast }(T(x),T(y))h(T(x),T(y))\leq \varphi (d(x,y))$\textit{\ for
all }$x,y\in X$\textit{, where }$h$\textit{\ stands for Hausdorff-Pompeiu
metric.}

\bigskip

\textbf{Definition 5.3. }\textit{For} $T:X\rightarrow \mathcal{P}_{cl}(X)$%
\textit{, where }$(X,d)$\textit{\ is a }$b$\textit{-metric space, we say
that }$T$ \textit{is a multivalued weakly Picard operator if for every }$%
x\in X$ \textit{and every }$y\in T(x)$\textit{\ there exists a sequence }$%
(x_{n})_{n\in \mathbb{N}}$\textit{\ of elements from }$X$\textit{\ such that:%
}

\textit{a) }$x_{0}=x$\textit{\ and }$x_{1}=y$\textit{;}

\textit{b) }$x_{n+1}\in T(x_{n})$\textit{\ for every }$n\in \mathbb{N}$%
\textit{;}

\textit{c) }$(x_{n})_{n\in \mathbb{N}}$\textit{\ is convergent and its limit
is a fixed point of }$T$\textit{.}

\bigskip

Now we can state our result.

\bigskip

\textbf{Theorem 5.1.} \textit{Let us consider} $T:X\rightarrow \mathcal{P}%
_{cl}(X)$, $\alpha :X\times X\rightarrow \lbrack 0,\infty )$ \textit{and }$%
\varphi \in \Gamma ^{\gamma }$\textit{, where }$(X,d)$\textit{\ is a
complete }$b$\textit{-metric space with constant }$s>1$\textit{\ and }$%
\gamma >\log _{2}s$\textit{, such that:}

\textit{a) }$T$\textit{\ is an }$\alpha _{\ast }$-$\varphi $\textit{%
-contractive multivalued operator of type }$(b)$;

\textit{b) }$T$ \textit{is }$\alpha _{\ast }$\textit{-admissible;}

\textit{c) there exists }$x_{0}\in X$\textit{\ and }$x_{1}\in T(x_{0})$%
\textit{\ satisfying the inequality }$\alpha (x_{0},x_{1})\geq 1$\textit{;}

\textit{d) for every convergent sequence} $(y_{n})_{n\in \mathbb{N}}$\textit{%
\ of elements from }$X$\textit{\ having limit }$y$\textit{, the following
implication is valid: }$(\alpha (y_{n},y_{n+1})\geq 1$ \textit{for every} $%
n\in \mathbb{N})\Rightarrow (\alpha (y_{n},y)\geq 1$ \textit{for every} $%
n\in \mathbb{N})$.

\textit{Then }$T$\textit{\ has a fixed point.}

\textit{Proof}. The same line of arguments given in the proof of Theorem 1
from [9] gives us a sequence $(x_{n})_{n\in \mathbb{N}}$\textit{\ }of
elements from $X$, with $x_{0}\neq x_{1}$, such that:

i)\textit{\ }$x_{n+1}\in T(x_{n})$\textit{\ }for every\textit{\ }$n\in 
\mathbb{N}$;

ii)\textit{\ }$\alpha (x_{n},x_{n+1})\geq 1$ for every $n\in \mathbb{N}$;

iii)\textit{\ }$d(x_{n},x_{n+1})\leq \varphi ^{\lbrack n]}(d(x_{0},x_{1}))$
for every $n\in \mathbb{N}$.

As $\varphi \in \Gamma ^{\gamma }$, the series $\overset{\infty }{\underset{%
n=1}{\sum }}n^{\gamma }\varphi ^{\lbrack n]}(d(x_{0},x_{1}))$ is convergent,
so, taking into account the comparison test and iii), we came to the
conclusion that $\overset{\infty }{\underset{n=1}{\sum }}n^{\gamma
}d(x_{n},x_{n+1})$ is convergent. Consequently, according to Lemma 2.2, $%
(x_{n})_{n\in \mathbb{N}}$ is Cauchy$.$

The same arguments as the ones used in the proof of Theorem 1 from [9]
assure us that the limit of $(x_{n})_{n\in \mathbb{N}}$ is a fixed point of $%
T$. $\square $

\bigskip

\textbf{Remark 5.1.} \textit{If hypothesis c) is replace by the following
condition: }$\alpha (x_{0},x_{1})\geq 1$\textit{\ for every }$x_{0}\in X$%
\textit{\ and every }$x_{1}\in T(x_{0})$\textit{, then the conclusion of the
above result is that }$T$\textit{\ is a multivalued weakly Picard operator.}

\bigskip

\textbf{Remark 5.2.}\textit{\ Using the same technique, one can generalize
Theorems 2 and 3 from }[9]\textit{.}

\bigskip

\textbf{References}

\bigskip

[1] A. Aghajani, M. Abbas and J.R. Roshan, Common fixed point of generalized
weak contractive mappings in partially ordered $b$-metric spaces, Math.
Slovaca, \textbf{64} (2014), 941-960.

[2] H. Aydi, M.F. Bota, E. Karapinar and S. Mitrovi\'{c}, A fixed point
theorem for set-valued quasi-contractions in $b$-metric spaces, Fixed Point
Theory Appl., 2012, 2012:88.

[3] I. A. Bakhtin, The contraction mapping principle in quasimetric spaces,
Funct. Anal., Unianowsk Gos. Ped. Inst., \textbf{30} (1989), 26-37.

[4] V. Berinde, Generalized contractions in quasimetric spaces, Seminar on
Fixed Point Theory, 1993, 3-9.

[5] V. Berinde, Sequences of operators and fixed points in quasimetric
spaces, Stud. Univ. "Babe\c{s}-Bolyai", Math., \textbf{16} (1996), 23-27.

[6] M. Boriceanu, M. Bota and A. Petru\c{s}el, Multivalued fractals in $b$%
-metric spaces, Cent. Eur. J. Math., \textbf{8} (2010), 367-377.

[7] M. Boriceanu, A. Petru\c{s}el and A.I. Rus, Fixed point theorems for
some multivalued generalized contraction in $b$-metric spaces, Int. J. Math.
Stat., \textbf{6} (2010), 65-76.

[8] M. Bota, A. Moln\'{a}r and C. Varga, On Ekeland's variational principle
in $b$-metric spaces, Fixed Point Theory, \textbf{12} (2011), 21-28.

[9] M. Bota,V. Ilea, E. Karapinar, O. Mle\c{s}ni\c{t}e, On $\alpha _{\ast
}-\varphi $-contractive multi-valued operators in $b$-metric spaces and
applications, Applied Mathematics \& Information Sciences, \textbf{9}
(2015), 2611-2620.

[10] D.W. Boyd and J.S. Wong, On nonlinear contractions, Proc. Amer Math.
Soc., \textbf{20} (1969), 458-464.

[11] F.E. Browder, On the convergence of succesive approximations for
nonlinear functional equations, Nederl. Akad. Wet., Proc., Ser. A 71, 
\textbf{30} (1968), 27-35.

[12] J. Caristi, Fixed point theorems for mappings satisfying inwardness
conditions, Tran. Amer. Math. Soc., \textbf{215} (1976), 241-251.

[13] C. Chifu and G. Petru\c{s}el, Fixed points for multivalued contractions
in $b$-metric spaces with applications to fractals, Taiwanese J. Math., 
\textbf{18} (2014), 1365-1375.

[14] S. Czerwik, Contraction mappings in $b$-metric spaces, Acta Math.
Inform. Univ. Ostraviensis, \textbf{1} (1993), 5-11.

[15] S. Czerwik, Nonlinear set-valued contraction mappings in $b$-metric
spaces, Atti Sem. Mat. Fis. Univ. Modena, \textbf{46} (1998), 263-276.

[16] J. Jachymski, Caristi's fixed point theorem and selections of
set-valued contractions, J. Math. Anal. Appl., \textbf{227} (1998), 55-67.

[17] M.A. Khamsi and N. Hussain, KKM mappings in metric type spaces,
Nonlinear Anal., \textbf{73} (2010), 3123-3129.

[18] W. A. Kirk, Caristi's fixed point theorem and metric convexity. Colloq.
Math., \textbf{36} (1976), 81-86.

[19] J. Matkowski, Integrable solutions of functional equations,
Dissertations Math. (Rozprawy), 127 (1976).

[20] S. K. Mohanta, Some fixed point theorems using $wt$-distance in $b$%
-metric spaces, Fasc. Math., \textbf{54} (2015), 125-140.

[21] H.N. Nashine and Z. Kadelburg, Cyclic generalized $\varphi $%
-contractions in $b$-metric spaces and an application to integral equations,
Filomat, \textbf{28} (2014), 2047-2057.

[22] M.O. Olatinwo, A fixed point theorem for multi-valued weakly Picard
operators in $b$-metric spaces, Demonstratio Math., \textbf{42} (2009),
599-606.

[23] M. P\u{a}curar, Sequences of almost contractions and fixed points in $b$%
-metric spaces, An. Univ. Vest Timi\c{s}. Ser. Mat.-Inform., \textbf{48}
(2010), 125-137.

[24] J.R. Roshan, N. Hussain, S. Sedghi and N. Shobkolaei, Suzuki-type fixed
point results in $b$-metric spaces, Math. Sci. (Springer), \textbf{9}
(2015), 153--160.

[25] B.D. Rouhani and S. Moradi, Common fixed point of multivalued
generalized $\varphi $-weak contractive mappings, Fixed Point Theory Appl.,
2010, Article ID 708984.

[26] I. A. Rus, Generalized $\varphi $-contractions, Math., Rev. Anal. Num%
\'{e}r. Th\'{e}or. Approximation, Math., \textbf{47} (1982), 175-178.

[27] S. Shukla, Partial $b$-metric spaces and fixed point theorems,
Mediterr. J. Math., \textbf{11} (2014), 703-711.

[28] H. Yingtaweesittikul, Suzuki type fixed point for generalized
multi-valued mappings in $b$-metric spaces, Fixed Point Theory and
Applications, 2013, 2013:215.

\bigskip

University of Bucharest

Faculty of Mathematics and Computer Science

Str. Academiei\ 14, 010014 Bucharest, Romania

E-mail: miculesc@yahoo.com, mihail\_alex@yahoo.com

\end{document}